\documentclass{roadef}

\usepackage{amsmath}

\begin{document}


\title{Improving Clique Decompositions of Semidefinite Relaxations for Optimal Power Flow Problems}

\def\shorttitle{Titre court}

\author{Julie Sliwak\inst{1}\inst{2}\inst{3}, Miguel F. Anjos\inst{4}, Lucas L\'etocart\inst{2} , Jean Maeght\inst{1}, Emiliano Traversi\inst{2}}

\institute{
RTE, Paris La Défense, France \\
\email{firstname.lastname@rte-france.com}
\and
LIPN, UMR CNRS 7030 - Université Paris 13, Villetaneuse, France \\

\email{\{letocart, traversi\}@lipn.univ-paris13.fr}
\and
Polytechnique Montréal, Département de mathématiques et de génie industriel, Montréal, Canada \\
\and
University of Edinburgh, School of Mathematics, Edinburgh, United Kingdom
\email{anjos@stanfordalumni.org}
}

\maketitle
\thispagestyle{empty}

\begin{abstract}
SSemidefinite Programming (SDP) provides tight lower bounds for Optimal Power Flow problems. However, solving large-scale SDP problems requires exploiting sparsity. In this paper, we experiment several clique decomposition algorithms that lead to different reformulations and we show that the resolution is highly sensitive to the clique decomposition procedure. Our main contribution is to demonstrate that minimizing the number of additional edges in the chordal extension is not always appropriate to get a good clique decomposition. 
\end{abstract}


\newpage

\section{Introduction}

The Alternating Current Optimal Power Flow (ACOPF) problem is a nonconvex optimization problem whose purpose is to provide an operating point of the power network minimizing generation costs. Solving this problem is essential for RTE, the French transmission system operator, because ACOPF problems are subproblems of decision problems involving millions of euros, e.g. operational planning and grid development. Yet there is still no efficient method to solve ACOPF problems to global optimality.  

Many convex relaxations have been proposed to evaluate the quality of feasible solutions computed with local methods  \cite{molzahn_survey_OPF}. SDP relaxations often provide tight lower bounds \cite{Zero_duality_gap} that are useful to prove global optimality. Some promising global optimization methods are based on these relaxations: Godard et al. propose an adaptation of the Mixed-Integer Quadratic Convex Reformulation method to ACOPF problems in \cite{B&BHadrien}, Gopalakrishnan et al. present a branch-and-bound approach using SDP relaxations in \cite{OPF_B&B_SDP} and Josz et al. apply the Lasserre hierarchy in \cite{josz_SOS_OPF} to achieve global optimality. All these methods depend on large-scale SDP problems being solved efficiently.

Many algorithms are proposed in the literature to solve SDP problems. In this paper, we focus on interior-point methods that are the most reliable and accurate algorithms as far as we know. Interior-point methods solve efficiently small-to-medium-sized SDP problems but do not scale well for large-scale problems because of the Hessian equation that implies forming and factorizing a fully-dense matrix at each step. Yet, sparse large-scale problems can be tackled exploiting sparsity. In this paper, we focus on clique decomposition techniques \cite{fukudaI, fukudaII, chordal_graphs_SDP} that perform well for ACOPF problems. Jabr \cite{jabr2012sparsity} and Molzahn et al. \cite{Molzahn_SDP} have each proposed a way to use clique decompositions for ACOPF problems and both propose to work on problems formulated in complex variables. However, there are many other possible clique decomposition procedures that provide different reformulations of a SDP problem. Most of the clique decomposition algorithms seek to minimize the number of added edges in the chordal extension. It results in decompositions with many small cliques for ACOPF problems and many linking constraints are required to handle overlaps between cliques. As linking constraints can slow down the resolution, Nakata et al. \cite{fukudaII} and Molzahn et al. \cite{Molzahn_SDP} propose heuristics for merging cliques in order to reduce their number. Yet, merging cliques means adding more edges in the chordal extension. Therefore, classical clique decomposition approaches may not always be the most appropriate since they focus on the number of added edges, regardless of important criteria such as the size of the largest clique or the number of linking constraints.

In this paper, we show that different clique decompositions are not equivalent in terms of resolution by comparing different chordal extension algorithms on RTE \cite{ACOPFdataRTE} and other MATPOWER \cite{matpower} datasets. Our main contribution is to demonstrate in two ways that computing a chordal extension with a minimum number of added edges is not a relevant choice for ACOPF problems. The first test consists in a comparison between reformulations coming from the SDP problem formulated in complex variables and reformulations coming from the same problem but formulated in real variables. The conclusions of this comparison are that less edges are necessary in the real case but the reformulations are more performant in the complex case. The second test relies on a new clique combination algorithm that improves a given decomposition by adding more edges in the chordal extension. The edges are added in such a way as to minimize the number of linking constraints while keeping small cliques. Computational tests show that this algorithm significantly speeds up the resolution of the largest OPF instances. Both tests highlight the fact that the number of added edges is not always the best criterion to minimize, which proves that the computation of the chordal extension is a question that deserves to be investigated.

This paper is organized as follows. Section 2 describes briefly the OPF problem and the classical rank relaxation. Common clique decomposition techniques are presented in section 3 along with computational comparisons. Section 4 compares reformulations coming from the complex SDP problem and reformulations coming from the real SDP problem. Section 5 details our clique combination algorithm and resulting enhancements. Section 6 concludes the paper.

\section{ACOPF problem and rank relaxation}

In this section, we present briefly the ACOPF formulation and its classical semidefinite rank relaxation.

\subsection{ACOPF formulation}

The ACOPF problem is defined on a power transmission network that can be modelled as an oriented graph $T(N,B)$ where $N$ represents the electrical buses and $B$ the branches (transmission lines, transformers). Let us denote $G \subset N$ the subset of generator buses. The ACOPF problem is defined as follows:

\begin{equation}\label{ACOPF}\tag{ACOPF}
\begin{aligned}
{\text{min }} & \sum_{g\in G} \mathbf{c_g}Re(S_g^{gen})+\mathbf{k_g}\\
\text{s.t. }
 & S^{gen}_n = \mathbf{S^{load}_n} + \sum_{b\in B^-(n)}S^{dest}_b(v)\\
 & +\sum_{b\in B^+(n)}S^{orig}_b(v) & \forall n\in N\\
  &(\mathbf{v_n^{min}})^2\leq |v_n|^2\leq (\mathbf{v_n^{max}})^2 & \forall n\in N\\
  & \mathbf{P_g^{min}}\leq Re(S_g^{gen})\leq \mathbf{P_g^{max}} & \forall g\in G\\
  & \mathbf{Q_g^{min}}\leq Im(S_g^{gen})\leq \mathbf{Q_g^{max}} & \forall g\in G\\
  & |i^{orig}_b(v)|^2\leq (\mathbf{i_b^{max}})^2 & \forall b\in B\\
  & |i^{dest}_b(v)|^2\leq (\mathbf{i_b^{max}})^2 & \forall b\in B\\
  & S^{gen}_n=0 & \forall n \in N\backslash G\\
& v_n \in \mathbb{C},\ S^{gen}_n \in \mathbb{C} & \forall n \in N\\
\end{aligned}
\end{equation} 

where $Re(z)$ stands for the real part of the complex number $z$ and $Im(z)$ for the imaginary part. All constant parameters are in bold. $\mathbf{c_g}$ and $\mathbf{k_g}$ represent respectively the linear cost and the constant cost of the generator bus $g \in G$. $\mathbf{S^{load}_n}$ represents the load. $\mathbf{v_n^{min}}$ and $\mathbf{v_n^{max}}$ are bounds on the voltage magnitude at bus $n$. $\mathbf{P_g^{min}}$ and $\mathbf{P_g^{max}}$  (respectively $\mathbf{Q_g^{min}}$ and $\mathbf{Q_g^{max}}$) are bounds on the active (respectively reactive) power at generator bus $g \in G$. $\mathbf{i_b^{max}}$ represents the current limit for branch $b\in B$. The variables are the voltage $v_n$ and the power $S^{gen}_n$ at each bus $n \in N$ with $S^{gen}_n=0$ for all buses $n \in N\backslash G$. Both are complex variables.

$B^-(n)$ is the set of entering branches at bus $n$ and $B^+(n)$ the set of exiting branches. $o(b)$ stands for the origin of branch $b$, $d(b)$ for the destination of branch $b$. The currents $i^{orig}_b(v)$ and $i^{dest}_b(v)$ of a branch $b\in B$ are linear functions of $v_{o(d)}$ and $v_{d(b)}$ depending on physical characteristics. The power on the lines are defined as follows:


\begin{equation}
\begin{array}{l}
S^{orig}_b(v)=v_{o(b)}\overline{i^{orig}_b(v)}\\[2mm]
S^{dest}_b(v)=v_{d(b)}\overline{i^{dest}_b(v)}
\end{array}
\label{S}
\end{equation}

%

Note that we use thermal limits modelled with current and linear generation costs to have a formulation closer to what is used in practice at RTE. 


\subsection{Rank relaxation}

Since costs are linear, the ACOPF problem can be written in terms of voltage variables only. The resulting problem is a Quadratically Constrained Quadratic Program (QCQP) with complex variables that can be expressed in a compact way:

\begin{equation}
\left\{
    \begin{array}{ll}
        min &v^HQ_0v \\
        s.t.& v^HQ_pv\leq a_p\ \forall p=1..m\\
        	& v \in \mathbb{C}^n
    \end{array}
\right.
\label{QCQP}
\end{equation}

The classical rank relaxation is constructed by introducing the Hermitian matrix $W=vv^H$. This constraint is equivalent to $W\succeq0$ and $rank(W)=1$. The rank relaxation is the SDP problem obtained by suppressing the rank constraint:

\begin{equation}
\left\{
    \begin{array}{ll}
        min &Q_0\cdot W \\
        s.t.& Q_p\cdot W\leq a_p\ \forall p=1..m\\
        &W \succeq 0
    \end{array}
\right.
\label{SDP}
\end{equation}

This SDP problem with complex numbers can be converted into a SDP problem with real numbers using the rectangular representation and a symmetric matrix of size $2n \times 2n$.

\section{Clique decomposition analysis}

In this section, we review the clique decomposition techniques based on matrix completion introduced by Fukuda et al. in \cite{fukudaI, fukudaII}. Then we show that the resolution with interior-point methods is highly sensitive to the clique decomposition.

\subsection{Clique decomposition framework}

The clique decomposition relies on the matrix completion theorem \cite{chordal_graphs_SDP}. This allows to replace the Positive Semidefinite (PSD) constraint on a big matrix by several PSD constraints on smaller submatrices at the price of adding linking constraints between these submatrices. 

Let us define $A$ as the aggregate sparsity pattern of a SDP problem, i.e., the matrix of nonzeros entries in the data matrices (objective and constraint matrices). Let us denote $G^A$ as the graph associated to A. The matrix completion theorem can be applied if and only if $G^A$ is chordal. A chordal graph is a graph with no induced cycle of four vertices or more. A clique is a subset of vertices that are all connected together. A maximal clique is a clique that is not included in another clique. 

The general framework for clique decomposition contains three steps. The first step consists in computing a chordal extension H of the aggregate sparsity pattern $G^A$. The second step consists in determining the list of maximal cliques in H, denoted by $L = \{C_1,...,C_r\}$. These maximal cliques define the submatrices that must be positive semidefinite in the reformulated SDP problem. From the list of maximal cliques $L$, a clique graph W can be defined  with nodes corresponding to maximal cliques and weighted edges between each pair of vertices defined by the number of shared nodes between the two corresponding cliques. The third step consists in computing a clique tree U by computing a maximum-weight spanning tree for W. This clique tree allows to specify linking constraints between submatrices.

It is clear that the clique decomposition depends on the chordal extension H but it is not completely clear how to precisely define a "good" decomposition from a practical point of view. One reasonable choice is to use decompositions that minimize the number of additional edges. However, finding the minimal chordal extension is NP-complete \cite{chordalext-NP-complete}. Bergman et al. have recently proposed an exact model with exponentially many constraints  to find a minimal extension \cite{bergman2019minimum} but it is only usable for small instances. For this reason, several efficient heuristics for adding few edges are available in the literature \cite{survey_min_triangulations}. The most commonly used is based on a Cholesky factorization of $G^A$ whose rows and columns have been permuted according to a minimum or approximate minimum degree ordering \cite{AMD}. The ordering has a significant impact on the fill-in, that is, on the number of edges added to $G^A$. Once H is computed, the list of maximal cliques can be computed in linear time \cite{tarjan1984simple}. Finally, the most widely used exact algorithm to compute the clique tree is Prim's algorithm \cite{prim} but there are other algorithms and they can lead to different optimal clique trees. However, the decompositions obtained usually do not differ much: the number of linking constraints do not differ significantly from one tree to another because the trees have the same overall weight.

We show in the following subsection that different chordal extension algorithms give different decompositions and that these different decompositions are not equivalent as regards resolution time.

\subsection{Numerical comparison of two chordal extension heuristics for SDP relaxations of ACOPF formulated in complex variables}

There are several algorithms to compute chordal extensions and the clique decompositions can vary greatly from one algorithm to another, which has an impact on the resolution. In this subsection, we compare clique decompositions coming from two different algorithms: a Cholesky factorization with an AMD ordering and the Minimum Degree (MD) heuristic \cite{MDheuristic} that is based on a dynamical computation of a minimum degree ordering of the vertices. We have picked these two algorithms because they are among the best heuristics known. Both algorithms are applied on the SDP relaxations of ACOPF formulated in complex variables, therefore providing maximal cliques with complex variables. In order to solve the SDP problems, these cliques are finally converted into cliques with real variables by doubling them.

Table 1 presents some results of this comparison on MATPOWER instances with more than 1000 buses. The tests were carried out on a Processor Intel® Core™ i7-6820HQ CPU @2.70GHz using the module MathProgComplex.jl \cite{MathProgComplex_powertech} with Julia 1.0.3. and the SDP solver MOSEK 8.1.0.72 was used for the resolution. We also used the packages JuMP.jl \cite{JuMP}, LightGraphs.jl \cite{LightGraphs.jl} and MetaGraphs.jl.

Let $nc^{MD}$, $nlc^{MD}$ and $t^{MD}$ be respectively the number of maximal cliques, the number of linking constraints and MOSEK resolution time obtained with the MD heuristic. Similarly, let $nc^{AMD}$, $nlc^{AMD}$ and $t^{AMD}$ be respectively the number of maximal cliques, the number of linking constraints and MOSEK resolution time obtained with Cholesky and AMD. Table~\ref{different_orderings} presents some ratios between these three quantities. All results are presented in percentages.

\begin{table}[!ht]
\centering

\label{different_orderings}
\begin{tabular}{|c|c|c|c|}
\hline
\rule[-7pt]{0pt}{20pt}
Instance  & $\frac{nc^{MD}}{nc^{AMD}}-1$  & $\frac{nlc^{MD}}{nlc^{AMD}}-1$ & $\frac{t^{MD}}{t^{AMD}}-1$\\[2mm]
\hline
case1888rte & 0.17\% & 1.7\% & 6.1\%\\
case3012wp &  -0.14\%& -3.3\% & -15.1\%\\
case6468rte &  0.10\%&0.31\% & -6.8\%\\
case9241pegase & -0.02\%& 4.9\% & 62\%\\
case13659pegase & -0.04\%& 4.7\% & 48\%\\
\hline 
\end{tabular}
\caption{Comparison of clique decompositions for MD and AMD }
\end{table}

This table is not exhaustive but some conclusions can be drawn from it. First, clique decompositions are different from one algorithm to another. In particular, there are large differences in resolution time even if the decompositions have about the same number of cliques. Theses differences can be explained by the differences in the number of linking constraints, especially for the biggest instances. Another factor can be the size of the biggest clique in the decomposition. It is also possible to generate decompositions which differ both in the number of cliques and in the number of linking constraints so differences in resolution time can be much more impressive. In particular, in preliminary tests we experienced that there exists "bad" clique decompositions by computing several "bad" orderings for Cholesky (e.g. random ordering, maximum degree ordering). In this case, a "bad" clique decomposition means a decomposition with many linking constraints, which can lead to memory issues with interior-point methods. This confirms the idea that the computation of the chordal extension is a point worth studying.

The results shown in this section are a good motivation
to continue exploring the impact of the clique reformulation on the resolution. The next section focus on the impact of computing the clique decomposition on SDP relaxations of ACOPF problems formulated in complex variables or on SDP relaxations of ACOPF problems formulated in real variables.

\section{Comparison of clique decompositions computed in the complex and in the real framework}
SDP relaxations of ACOPF problems are naturally expressed with complex numbers but in practice, in order to be solved by the available solvers, they are rewritted in terms of only real variables. The process of rewriting the model with real variables leads to a different (and bigger) aggregate sparsity pattern A. Such a A could be itself decomposed in a way that would not be possible with the original formulation
with the complex variables. Therefore, there are two possibilities for the clique decomposition: either applying the procedure on the complex SDP formulation and converting the complex cliques to real cliques or directly applying the clique decomposition procedure on the real SDP formulation. However, computing the chordal extension on the complex or the real problem is not equivalent. We first show it theoretically on a small MATPOWER \cite{matpower} example, LMBM3. 

Let us define $G^c$ (respectively $G^r$) as the graph associated to the aggregate sparsity pattern A of the complex (respectively real) SDP relaxation. For any ACOPF instance as defined in section 2, $G^c$ is the network graph $T(N,B)$. Let us denote $H^c$ (respectively $H^r$) the chordal extension computed from $G^c$ (respectively $G^r$) with a given algorithm. The chordal extension $H^c$ can be converted to real numbers in the same way as the complex graph $G^c$ can be converted to the real graph $G^r$. Let us denote $H^c_{real}$ this conversion of $H^c$. $H^c_{real}$ is also a chordal extension of $G^r$. The objective of this section is to compare $H^c_{real}$ and $H^r$.

For LMBM3, $G^c$ is a clique of size 3, denoted by $K_3$ and $G^r$ is the graph with the 6 following nodes $V^r= \{1_{Re}, 1_{Im},2_{Re}, 2_{Im},3_{Re}, 3_{Im} \}$ and the 12 following edges: 
\begin{equation*}
\begin{array}{ll}
 E^r = & \{(1_{Re},2_{Re}), (1_{Re},2_{Im}), (1_{Im},2_{Re}), (1_{Im},2_{Im}),\\ 
 & (1_{Re},3_{Re}), (1_{Re},3_{Im}), (1_{Im},3_{Re}), (1_{Im},3_{Im}),\\
 & (2_{Re},3_{Re}), (2_{Re},3_{Im}), (2_{Im},3_{Re}), (2_{Im},3_{Im}) \}.
\end{array}
\end{equation*}

Since $G^c$ is complete and thus chordal, $H^c=G^c=K_3$ regardless of the chordal extension algorithm used. Clearly, there is one maximal clique with 3 nodes and no linking constraints. When converted to real numbers, we get $H^c_{real} = K_6$, the clique of size 6, which accounts to adding 3 edges to $G^r$. This chordal extension gives one maximal clique with 6 real nodes and no linking constraints. 
On the other hand, $G^r$ is not chordal because there are several induced subgraphs of size 4, e.g. $1_{Re} - 2_{Im} - 1_{Im} - 2_{Re}$. It suffices to add the two edges $(2_{Re},2_{Im})$ and $(3_{Re},3_{Im})$ to $G^r$ to get a minimal chordal extension $H^r$ \footnote{There are other possibilities for $H^r$, e.g., if the edges $(1_{Re},2_{1m})$ and $(3_{Re},3_{Im})$ are added, but the results will be the same.}.  $H^r$ is chordal because $1_{Re} - 1_{Im} - 2_{Re} - 2_{Im} - 3_{Re} - 3_{Im}$ is a perfect elimination ordering. This chordal extension gives two maximal cliques: $\{2_{Re},2_{Im}, 3_{Re}, 3_{Im}, 1_{Re}\}$ and $\{2_{Re},2_{Im}, 3_{Re}, 3_{Im}, 1_{Im}\}$ and 10 linking constraints.

This small example proves that it is not equivalent to compute the clique decomposition on the complex problem or on the real problem from the theoretical point of view: more edges are added in the complex case to get a chordal extension. This idea is confirmed by computational tests that show that there are on average 50\% more edges added when the chordal extension is computed in the complex case rather than in the real case for MATPOWER instances with more than 1000 buses. However, adding less edges in the chordal extension is not necessarily advantageous from the numerical resolution point of view. Indeed, adding less edges in the chordal extension can mean having more cliques and more linking constraints as shown in Table 2. This table presents comparisons between clique decomposition coming from the complex case (computed on $G^c$) and clique decomposition coming from the real case (computed on $G^r$). More precisely, both clique decompositions are based on a chordal extension obtained by a Cholesly factorization with an Approximate Minimum Degree (AMD) \cite{AMD} ordering. Let $L^c$ (respectively $L^r$) be the list of maximal cliques obtained from the chordal extension $H^c$ (respectively $H^r$). The second column represents the ratio for the number of cliques defined as $\frac{|L^r|}{|L^c|}-1$. Similarly, the third column represents the ratio for the number of linking constraints, i.e., $\frac{nlc^r}{nlc^c}-1$ with $nlc^c$ (respectively $nlc^r$) the number of linking constraints when the clique decomposition is done for the complex (respectively real) problem. Finally, the fourth column represents the ratio for MOSEK resolution time, i.e., $\frac{t^r}{t^c}-1$ with $t^c$ (respectively $t^r$) the resolution time with the clique decomposition coming from the complex (respectively real) problem. All results are presented in percentages and the tests were carried out in the same manner as in subsection 3.2.

\begin{table}[!ht]
\centering

\label{C_vs_R}
\begin{tabular}{|c|c|c|c|}
\hline
\rule[-7pt]{0pt}{20pt}
Instance  & $\frac{|L^r|}{|L^c|}-1$  & $\frac{nlc^r}{nlc^c}-1$ & $\frac{t^r}{t^c}-1$\\[2mm]
\hline
case1888rte & 51\% & 24\% & 17\%\\
case3012wp &  39\%& 15\% & 4\%\\
case6468rte &  40\%&19\% & 58\%\\
case9241pegase &  48\%& 18\% & 28\%\\
case13659pegase &  59\%& 20\% & 23\%\\
\hline 
\end{tabular}
\caption{Comparison of clique decompositions on $G^c$ or $G^r$}
\end{table}

This table shows that the clique decomposition procedure on $G^r$ gives approximatly 50\% more cliques and 20\% more linking constraints than on $G^c$, which slows down the resolution (the computing time increases by at least 4\% more and up to 58\% more). Therefore, even if less edges are added to the graph $G^r$ when working directly in real numbers, it seems better to work in complex numbers. This comparison demonstrates that computing a chordal extension with a minimum number of added edges is not necessarily profitable. The next section presents a clique combination algorithm that leads to the same conclusion.

\section{A new clique combination approach}

Several heuristics do already exist to compute chordal extensions  but to the best of our knowledge, they all try to minimize the number of added edges. However, nothing proves that a chordal extension with few added edges gives a decomposition easier to solve. In contrast, minimizing the size of the biggest clique or the number of linking constraints could be more interesting because these quantities have an impact on the resolution. To support this claim, we propose a new clique combination algorithm that allows to reduce the number of linking constraints in a given decomposition while keeping small matrices. This algorithm consists in merging some cliques that are adjacent in the clique tree by solving the problem (IP) below. There is one binay variable per edge in the clique tree. If this variable equals 1, it means that the cliques at the ends of the edge will be merged. We allow at most one combination per clique and combinations that result in a clique of size greater than a certain size are forbidden. The objective is to maximize the number of linking constraints that will be eliminated thanks to the combinations. More precisely, the model is the following:

\begin{equation}
    \begin{array}{ll}
        \max & 0.5\sum_{e\in E}  nc(C_e)(nc(C_e)+1)x_e\\[2mm]
       s.t.  & size(C_e)x_e \leq S^{max} \ \forall e\in E\\[2mm]
            & \sum_{e\in E(i)} x_e \leq 1 \ \forall i \in L \\[2mm]
            & x_e \in \{0,1\}\ \forall e\in E
    \end{array}
    \tag{IP}
\end{equation} with
\begin{itemize}
	\item $U=(L,E)$ the clique tree computed for a given decomposition ;
    \item $nc(C_e)= |C_i \bigcap C_j|$ the number of common nodes for $e=(C_i,C_j)\in E$;
    \item $size(C_e) = |C_i|+|C_j|-|C_i \bigcap C_j|$ the size of the merged clique for $e=(C_i, C_j) \in E$;
    \item $S^{max}$ a parameter that defines the maximal size of merged cliques;
    \item $E(i)$ the set of edges that have $i$ as extremity.
\end{itemize}

\begin{table}[h]
\begin{center}
\begin{tabular}{|c|c|c|c|c|c|c|c|}
\hline
Number of times (IP) is applied & 0& 1 & 2 & 3 & 4 & 5 & 6\\
\hline
Total resolution time (seconds) & 3397 & 2316 & 1807 &1691 &1533 & 1710 & 1826\\
\hline 
\end{tabular}
\end{center}
\caption{Total MOSEK resolution time for MATPOWER instances}
\end{table} 

This algorithm has been tested with the solver Xpress for $S^{max}=50$ and the results are presented in Table 3. This table represents the evolution of the total resolution time depending on the number of times the combination algorithme (IP) is used. This table shows that applying this algorithm at least once improves significatively MOSEK resolution time and the best time is obtained when the algorithm is applied 4 times. Therefore, it is worth combining cliques, i.e., adding more edges in the chordal extension, to improve the performance of a given decomposition for ACOPF problems.

\section{Conclusion and future work}
Sections 4 and 5 show that, for ACOPF problems, there is no real reason to build chordal extensions with as few added edges as possible since adding more edges can lead to better performance. These results encourage us to change the approach for computing chordal extensions and in the future we will focus on the formulation of an adapted heuristic that computes performant chordal extensions for ACOPF problems.

\bibliographystyle{plain}
\bibliography{D:/Users/sliwakjul/Documents/Biblio/biblio}

\end{document}